\documentclass[12pt]{article}

\usepackage{amsmath}
\usepackage{amssymb}
\usepackage{graphicx}
\usepackage{epsfig}

\usepackage{hyperref} 

\newtheorem{theorem}{Theorem}
\newtheorem{lemma}[theorem]{Lemma}

\def \S {{\cal{S}}}
\def \B {{\cal{B}}}
\def \A {{\cal{A}}}
\def \a {\alpha}
\def \b {\beta}

\def \deg {{\rm deg}}
\def \proof {\noindent{\bf Proof}\quad}

\def \eps {\varepsilon}

\def \leq {\leqslant}
\def \geq {\geqslant}

\def\endmark{\hskip 2em$\square$\par}
\def \qed {\hfill\endmark}

\title{\bf Hamilton path decompositions\\ of complete multipartite graphs\footnote{This research was supported under the Australian Research
Council's Discovery Projects funding scheme, project numbers
DP150100506 and DP150100530.}}

\author{
\hspace{-0.7cm}
\begin{minipage}[c]{15cm}
\begin{center}
\begin{tabular}{ccc}
Darryn Bryant & Hao Chuien Hang & Sarada Herke\\
db@maths.uq.edu.au & hanghc@hotmail.com & s.herke@uq.edu.au\\
&&\\
\end{tabular}
School of Mathematics and Physics \\ 
University of Queensland \\
Qld 4072, Australia\\
\end{center}
\end{minipage}
}

\date{ }

\begin{document}
\maketitle\thispagestyle{empty}
\def\baselinestretch{1.5}\small\normalsize

\begin{abstract}
We prove that a complete multipartite graph $K$ with $n>1$ vertices and $m$ edges 
can be decomposed into edge-disjoint Hamilton paths if and only if 
$\frac m{n-1}$ is an integer and the maximum degree of $K$ is at most $\frac {2m}{n-1}$. 
\end{abstract}

There has been interest in problems concerning decompositions of graphs into Hamilton cycles, or into Hamilton paths, for many years. 
A well-known construction of Walecki (see \cite{Als,Luc}) can be used to obtain a decomposition of any complete graph of odd order into Hamilton cycles, and a decomposition of any complete graph of even order into Hamilton paths. 
A complete multipartite graph has its vertices partitioned into parts and two vertices are adjacent if and only if they are from distinct parts. 
In 1976, Laskar and Auerbach \cite{LasAue} showed that a complete multipartite graph can be decomposed into Hamilton cycles if and only if
it is regular of even degree, a condition which is obviously necessary. 
In this paper we prove a corresponding result for decompositions of complete multipartite graphs into Hamilton paths.

If $K$ is any graph with $n>1$ vertices and $m$ edges that can be decomposed into Hamilton paths, 
then clearly $t=\frac m{n-1}$ is an integer,
equal to the number of Hamilton paths in the decomposition, 
and the degree of each vertex of $K$ is at most $2t$,
because each Hamilton path has maximum degree 2.
We say that a complete multipartite graph with $n>1$ vertices and $m$ edges is 
{\em admissible} if and only if it satisfies these conditions.
The trivial complete multipartite graph with one vertex is also admissible. 
Our main result is that a complete multipartite graph has a decomposition into Hamilton paths if and only if it is 
admissible (see Theorem \ref{mainresult}).

Admissible complete multipartite graphs are plentiful and arise in no obvious regular pattern, 
although there are a few easily described infinite families.  The number of admissible complete multipartite graphs with at least two parts and having order at most $20$, at most $40$, and at most $60$ is $53$, $275$, and $917$, respectively.  Below we list some examples of admissible complete multipartite graphs,
including a few infinite families.  
We use the notation $K_{A_1,\ldots,A_r}$, where $A_1,\ldots,A_r$ are non-empty sets,
to denote the complete multipartite graph with parts 
$A_1,\ldots,A_r$, and we use the notation 
$K_{a_1^{x_1}, a_2^{x_2}, \ldots, a_s^{x_s}}$, where $a_1,\ldots,a_s$ and $x_1,\ldots,x_s$ are positive integers, to denote any complete multipartite graph having $x_i$ parts of cardinality $a_i$ for $i=1,\ldots,s$.
In the list below we write just $a_1^{x_1}, a_2^{x_2}, \ldots, a_s^{x_s}$ rather than $K_{a_1^{x_1}, a_2^{x_2}, \ldots, a_s^{x_s}}$.

\vspace{0.3cm}

\noindent
$
\begin{array}{lllll}
1^4,2,3; & 2^3,3^2,4; & 3^5,4,6;  & 6^6,7^2,10; & 8^6,10^4;\\
9^8,10^2,12,13;\quad & 10^6,11^{27},12^2,13;\quad & 9^{38},10^2,12,13^2,21;\quad & 57^{24},61^2,63^3;\quad & 75^{19},76^{22},78,79,87;\\
\end{array}
$
$
\begin{array}{lllll}
1^x		\text{ where } x\geq 2 \text{ is even;}&&&
a, (a+1)^x  \text{ where } a,x\geq 1\text{ and }(a+1)x \text{ is even;}& \\
1^{(a-1)^2}, a \text{ for } a \geq 2; \hspace{2.5cm} &&&
a^{2a-1}, 2a  \text{ for } a\geq 1.\qquad &\\
\end{array}
$

\vspace{0.5cm}

Throughout the paper, our graphs are allowed to have edges of multiplicity greater than $1$ and loops are permitted. 
A graph with no edges of multiplicity greater than $1$ and with no loops is called {\em simple}. Let $G$ be a graph. The number of vertices,
number of edges and number of components in $G$ are denoted by $v(G)$, $e(G)$ and $c(G)$, respectively, and the degree in $G$ of the vertex $x$ is denoted by $\deg_G(x)$. 
We assume that a loop on vertex $x$ contributes $2$ to the degree of $x$.  
The number of edges in a path is called its {\em length}.

A decomposition of a graph $K$ is a set $\{G_1,\ldots,G_t\}$ of subgraphs of $K$ such that $E(G_i)\cap E(G_j)=\emptyset$ for 
$1\leq i<j\leq t$ and $E(G_1)\cup\cdots\cup E(G_t)=E(K)$. 
A spanning subgraph is called a {\em factor} and a decomposition into factors is called a {\em factorisation}. 
If $S\subseteq V(G)$, then $G$ is {\em almost regular on $S$} if $|\deg_G(x)-\deg_G(y)|\leq 1$ for all $x,y\in S$,
and $G$ itself is {\em almost regular} if $G$ is almost regular on $V(G)$. 
In an edge colouring of $G$ with colours $c_1,\ldots,c_t$, 
the {\em colour class} $c_i$ refers to the spanning subgraph of $G$ whose edges are assigned colour $c_i$. 
An edge colouring is {\em almost regular} on $S$ if each of its colour classes is almost regular on $S$. 

Let $G$ be a graph. 
The {\em degree partition} of $V(G)$ is the partition where vertices $x$ and $y$ are in the same part if and only 
if they have the same degree. 
For any partition $\S=\{S_1,\ldots,S_r\}$ of $V(G)$, we define the graph $G^\S$ as follows. 
The vertex set of $G^\S$ is 
$\{S_1,\ldots,S_r\}$ and for each edge $e$ of $G$ there is a corresponding edge $e^\S$ in 
$G^\S$ such that the endpoints of $e^\S$ are $S_x$ and $S_y$ if and only if the endpoints of $e$ are in $S_x$ and $S_y$.

A brief outline of the proof of our main theorem is as follows. 
We begin (see Lemma \ref{partsizes}) by constructing a factorisation of any admissible $K_{A_1,\ldots,A_r}$ with the property that each factor $F$ is the vertex disjoint union of cycles and a single path and $F^\B$ is connected, where $\B$ is the degree partition of the vertices of
$K_{A_1,\ldots,A_r}$. Thus, although each $F$ is not necessarily connected, if we group the vertices of $F$ according to the cardinalities of their respective parts in $K_{A_1,\ldots,A_r}$ 
(equivalently according to their degree in $K_{A_1,\ldots,A_r}$) and collapse each of these groups into a single vertex, then the resulting graph is connected. 
To achieve this, we prove that $K_{A_1,\ldots,A_r}^\B$ contains a spanning, or near-spanning, star of multiplicity 
$2t$ where $t$ is the number of Hamilton paths we require in our decomposition (see Lemma \ref{starlemma}), 
and the result is then an easy consequence of Petersen's $2$-factor Theorem and Lemma \ref{equitabilise}. 
A {\em star of multiplicity $\mu$}
is a graph with one central vertex $u$ and a (possibly empty) set $V$ of other vertices
such that $u$ is joined to each vertex in $V$ with an edge of multiplicity $\mu$. 

From the factorisation described in the preceding paragraph, we then use Lemma \ref{connectinglemma} to obtain a new factorisation of $K_{A_1,\ldots,A_r}$ with the property that each factor $F$ is the vertex disjoint union of cycles and a single path and $F^\A$ is connected, where $\A$ is the partition $\{A_1,\ldots,A_r\}$ (see Lemma \ref{partsgraph}). 
Thus, if we now group the vertices of $F$ according to their respective parts in 
$K_{A_1,\ldots,A_r}$ and collapse each of these groups into a single vertex, then the resulting graph is connected (so each factor is connected under a finer partition than previously).  
Finally, in the proof of Theorem \ref{mainresult} 
we again use Lemma \ref{connectinglemma} to modify our factorisation so that each
factor becomes connected. Thus, each factor consists of a single path and hence is a Hamilton path.

We now briefly discuss some existing work on path decompositions and some related results.
There is a well-known conjecture of Gallai from 1966 that every simple connected graph of order $n$ can be decomposed into $\lceil \frac n2\rceil$ or fewer paths.  In 1968, Lov\'{a}sz \cite{Lov} proved that every simple connected graph of order $n$ can be decomposed into $\lfloor \frac n2\rfloor$ or fewer paths and cycles.  A consequence of Lov\'{a}sz's result is that Gallai's conjecture holds for graphs in which every vertex has odd degree.  See \cite{BotJim, HarMcG} and the references therein for some recent progress on Gallai's conjecture.

A lot is known about path decompositions of complete graphs. In 1983, Tarsi \cite{Tar} proved that a complete graph 
of order $n$ and multiplicity $\lambda$, denoted $\lambda K_n$, 
can be decomposed into paths of length $m$ if and only if $m\leq n-1$ and $m$ divides $e(\lambda K_n)$. 
Tarsi \cite{Tar} also conjectured that $\lambda K_n$ can be decomposed into $t$ paths of lengths 
$m_1, \dots, m_t$ if and only if $m_i\leq n-1$ for $i = 1, \dots, t$ and $m_1+\dots+m_t = e(\lambda K_n)$.  
Tarsi's conjecture was proved in 2010 \cite{Bry}.
Parker \cite{Par} has proved necessary and sufficient conditions 
for the existence of a decomposition of a complete bipartite graph into paths of any fixed length. 

By alternately colouring the edges of each path, from any 
decomposition of a graph $K$ into $t$ paths we can obtain 
a proper edge colouring of $K$ with $2t$ colours.  
We show in Lemma \ref{maxdeg2tlemma} that admissible complete multipartite graphs with $n$ vertices and $m$ edges are either complete graphs or have maximum degree 
$2t$ where $t=\frac m{n-1}$. Thus, by showing that admissible complete 
multipartite graphs of maximum degree $2t$ can be decomposed into $t$ 
Hamilton paths, we show that they have proper edge colourings with $2t$ colours. 
This represents a special case of the result of Hoffman and Rodger that 
complete multipartite graphs are {\em class 1} if and only if they are not 
{\em overfull} \cite{HofRod}.

Since Thomassen's paper \cite{Tho1}, which shows that 171-edge-connected graphs of size divisible by 3 can be decomposed into paths of length 3, 
there has been considerable interest in path decompositions of highly connected graphs
(also see \cite{Tho2}).
In \cite{BotMotOshWak}
it is shown that for any fixed length $k$, there is a constant $C(k)$ such that 
any $C(k)$-edge-connected graph $G$ has a decomposition into paths of length $k$ if and only if $k$ divides the size of $G$. 

We will use a theorem from \cite{Bry2} which concerns almost-regular edge colourings of hypergraphs. The following lemma is an immediate consequence of the result in \cite{Bry2}, as it applies to ordinary graphs. Also see \cite{BryMae} for the same result in the case of edge colourings of complete graphs. As noted in \cite{BryMae}, 
this result has some similarities and connections with the method of graph amalgamations \cite{AndRod}, and could be proved by that method.

\begin{lemma}\label{almostregularcolouring}{\rm {(\cite{Bry2})}}
If $G$ is a graph, $\gamma$ is an edge colouring of $G$, and $S\subseteq V(G)$ such that any permutation of $S$ is an automorphism of $G$, then there exists an edge colouring $\gamma'$ of $G$ such that 
\begin{itemize}
\item [$(a)$] for each colour $c$, the number of edges of colour $c$ is the same for $\gamma$ and $\gamma'$;
\item [$(b)$] for each $v\in V(G)\setminus S$ and each colour $c$, the number of edges of colour $c$ incident with $v$ is the same for $\gamma$ and $\gamma'$;
\item [$(c)$] the colour of any edge of $G$ that has both endpoints in $V(G)\setminus S$ is the same in $\gamma$ and $\gamma'$; and
\item [$(d)$] $\gamma'$ is almost regular on $S$.
\end{itemize}
\end{lemma}

The next two results are proved using Lemma \ref{almostregularcolouring}, and will be used in the proofs of Lemmas
\ref{partsizes} and \ref{partsgraph} and Theorem \ref{mainresult}. 

\begin{lemma}\label{equitabilise}
If $\{G_1,\ldots,G_t\}$ is any factorisation of $K_{A_1,\ldots,A_r}$ and $\A=\{A_1,\ldots,A_r\}$, 
then there exists a factorisation $\{F_1,\ldots,F_t\}$ of $K_{A_1,\ldots,A_r}$ such that for $i = 1, \dots, t$,
\begin{itemize}
\item [$(1)$] $F^\A_i=G^\A_i$; and
\item [$(2)$] for $j=1,\ldots,r$, $F_i$ is almost regular on $A_j$.
\end{itemize}
\end{lemma}

\proof
Define an edge colouring $\gamma$ of $K_{A_1,\ldots,A_r}$ with colours $c_1,\ldots,c_t$ by colouring the edges of $G_i$ with colour $c_i$ for 
$i=1,\ldots,t$. Now let $\gamma'$ be the new edge colouring of $K_{A_1,\ldots,A_r}$ obtained from $\gamma$ by applying Lemma \ref{almostregularcolouring}, first with $S=A_1$, then with $S=A_2$ and so on. 
Let $\{F_1,\ldots,F_t\}$ be the factorisation of $K_{A_1,\ldots,A_r}$ 
in which for $i=1,2,\ldots,t$,
$F_i$ contains precisely the edges that are assigned colour $c_i$ by $\gamma'$.
Property (1) follows immediately from (b) and (c) of Lemma \ref{almostregularcolouring} and  property (2) 
follows from (d) of Lemma \ref{almostregularcolouring}.
\qed

\begin{lemma}\label{connectinglemma}
Let $G$ be a graph, let $\a,\b\in V(G)$ such that the transposition $(\a\,\b)$ is an automorphism of $G$,
let $\{G_1,\ldots,G_t\}$ be a factorisation of $G$ such that 
for $i=1,\ldots,t$ either 
$\a$ and $\b$ are in the same component of $G_i$ or 
at least one of $\a$ and $\b$ is in a cycle of $G_i$.
Then there exists a factorisation $\{F_1,\ldots,F_t\}$ of $G$ such that for $i=1,\ldots,t$,
\begin{itemize}
\item [$(1)$] $e(F_i)=e(G_i)$; 
\item [$(2)$] for each $v\in V(G)\setminus\{\a,\b\}$, $\deg_{F_i}(v)=\deg_{G_i}(v)$; 
\item [$(3)$] $F_i-\{\a,\b\}=G_i-\{\a,\b\}$; 
\item [$(4)$] $F_i$ is almost regular on $\{\a,\b\}$; and 
\item [$(5)$] $\a$ and $\b$ are in the same component of $F_i$.
\end{itemize}
\end{lemma}

\proof
Define an edge colouring $\gamma$ of $G$ with colour set $\{c_1,\ldots,c_t\}\cup\{c'_1,\ldots,c'_t\}$ as follows.
For $i=1,\ldots,t$, 
\begin{itemize}
\item if $\a$ and $\b$ are in the same component of $G_i$, then the edges of a path from $\a$ to $\b$ in $G_i$
are assigned colour $c'_i$, and every other edge of $G_i$ is assigned colour $c_i$;
\item if $\a$ and $\b$ are in distinct components of $G_i$, then there is a cycle $C$ in $G_i$ 
that contains either $\a$ or $\b$, the edges of $C$ are assigned colour $c'_i$, and every other edge of 
$G_i$ is assigned colour $c_i$.
\end{itemize}
Let $\gamma'$ be the edge colouring of $G$ with colour set $\{c_1,\ldots,c_t\}\cup\{c'_1,\ldots,c'_t\}$ that is 
obtained from $\gamma$ by applying Lemma \ref{almostregularcolouring} with $S=\{\a,\b\}$ and for $i=1,\ldots,t$ let 
$F_i$ be the factor of $G$ containing the edges of colour $c_i$ or $c'_i$ in $\gamma'$.

It follows that the edges of colour $c'_i$ in $\gamma'$ form a path from $\a$ to $\b$ in $F_i$, and 
hence that $\a$ and $\b$ are in the same component of $F_i$. So we have proved that (5) holds. 
Properties (1), (2), (3) and (4), respectively, follow from (a), (b), (c) and (d), respectively, of Lemma \ref{almostregularcolouring}.
To see that (4) holds, observe that in $\gamma'$, there is exactly one edge of colour $c'_i$ incident with each of $\a$ and $\b$
and that colour class $c_i$ is almost regular on $\{\a,\b\}$. 
\qed

\vspace{0.3cm}

Lemmas \ref{maxdeg2tlemma} and \ref{epsboundlemma} are required only for the proof of Lemma \ref{starlemma}.  Lemma \ref{maxdeg2tlemma} was mentioned earlier, and implies that in any Hamilton path decomposition of a complete multipartite graph, that is not a complete graph, 
there is a vertex that is not the endpoint of any path. 
Lemma \ref{epsboundlemma} is a somewhat technical result that says that if there is a unique part of smallest cardinality $a$, then all the other parts have cardinality $a+1$, and that if there is more than one smallest part, then 
the number of vertices is at most $r(a+1)-2$ where $a$ is the cardinality of smallest part and $r$ is the number of parts.

\begin{lemma}\label{maxdeg2tlemma}
If $K_{A_1,\ldots,A_r}$ is admissible with $r\geq 2$, $n=v(K_{A_1,\ldots,A_r})$ and $t=e(K_{A_1,\ldots,A_r})/(n-1)$, 
then $K_{A_1,\ldots,A_r}$ is either a complete graph or it is not regular and has maximum degree $2t$.
\end{lemma}

\proof
Let $\Sigma$ denote the sum of the degrees of the vertices in $K_{A_1,\ldots,A_r}$. 
Since 
$\Sigma=2e(K_{A_1,\ldots,A_r})$ and $t=e(K_{A_1,\ldots,A_r})/(n-1)$ we have $2t=\frac\Sigma{n-1}$. 
If $K_{A_1,\ldots,A_r}$ is not a complete graph, then $\Sigma<n(n-1)$ from which it follows that 
$\Sigma>n(\frac\Sigma{n-1}-1)=n(2t-1)$. Thus, $K_{A_1,\ldots,A_r}$ has maximum degree $2t$. Also, 
having maximum degree $2t$ means that $K_{A_1,\ldots,A_r}$ is not regular, because if it were regular, then we would have 
$\Sigma=2tn$ which contradicts $2t=\frac\Sigma{n-1}$.
\qed

\vspace{0.3cm}

The argument used in the proof of Lemma \ref{maxdeg2tlemma} shows that every simple non-complete graph that has a decomposition into Hamilton paths is not regular, and that in any such decomposition there is a vertex that is not the endpoint of any path. 

\begin{lemma}\label{epsboundlemma}
Let $K_{A_1,\ldots,A_r}$ be admissible with $r\geq 2$, let $a_i=|A_i|$ for $i=1,\ldots,r$, and let 
$a_1\leq\cdots\leq a_r$.
If $a_1<a_2$, then $a_i=a_1+1$ for $i=2,\ldots,r$, and if $a_1=a_2$, then 
$a_1+\cdots+a_r\leq r(a_1+1)-2$.
\end{lemma}

\proof
If $K_{A_1,\ldots,A_r}$ is a complete graph, then the result holds (since $r\geq 2$) so we can assume otherwise.  
Let $a=a_1$ and let $p=|\{i:a_i=a, i=1,\ldots,r\}|$ and let $\sigma=\sum_{i=p+1}^r a_i$.
The number of vertices in $K_{A_1,\ldots,A_r}$ is $n=pa+\sigma$, the number of edges is 
$m={p\choose 2}a^2+pa\sigma+\sum_{p+1\leq i<j\leq r}a_ia_j$ and the maximum degree is 
$\Delta=(p-1)a+\sigma$.
Since $K_{A_1,\ldots,A_r}$ is admissible, we know that $t=\frac m{n-1}$ is an integer and by Lemma \ref{maxdeg2tlemma} we 
know that $\Delta=2t$ (since $K_{A_1,\ldots,A_r}$ is not a complete graph). Thus, $\Delta(n-1)=2m$. 
Substituting the above expressions for $n$, $m$ and $\Delta$ into this equation we obtain 
\[
((p-1)a+\sigma)(pa+\sigma-1)=p(p-1)a^2+2pa\sigma+2\sum_{p+1\leq i<j\leq r}a_ia_j
\]
which simplifies to 
\[\sum_{j = p+1}^r a_j^2 = \sigma + (p-1)a + a\sigma.\]
By the Cauchy-Schwarz Inequality, $\sigma^2 \leq (r-p) \sum\limits_{j = p+1}^{r} a_j^2 $, and so we have (using the fact that $r-p>0$,
which holds because Lemma \ref{maxdeg2tlemma} tells us that $K_{A_1,\ldots,A_r}$ is not regular)
\[
\frac{\sigma^2}{(r-p)} \leq \sigma + (p-1)a + a\sigma.
\]
which is equivalent to
\[
\sigma - (r-p)(a+1) \leq \frac{(p-1)(r-p)a}{\sigma}.
\]
If $a_1<a_2$, then $p=1$ and we have $\sigma\leq (r-1)(a+1)$ which implies $a_i=a_1+1$ for $i=2,\ldots,r$.
On the other hand, if $a_1=a_2$, then $p>1$, and using the fact that 
$(r-p)a < \sigma$ we obtain 
$
\sigma - (r-p)(a+1) < p-1
$
which means that 
$
\sigma - (r-p)(a+1)\leq p-2
$
and hence that  
$
\sigma+ap\leq r(a+1)-2.
$
The completes the proof because $\sigma+ap=a_1+\cdots+a_r$.
\qed

\begin{lemma}\label{starlemma}
Let $K=K_{A_1,\ldots,A_r}$ be admissible with $r\geq 2$, let $n=v(K)$, let $t=e(K)/(n-1)$,
and let $\B$ be the degree partition of $V(K)$.
Then $K^\B$ contains a star $S$ of multiplicity $2t$
such that $V(K^\B)\setminus V(S)$ is either empty or is $\{A_x\}$ for some $x\in\{1,\ldots,r\}$.
\end{lemma}

\proof
If $K$ is a complete graph, 
then $K^\B$ has a single vertex, and the result holds trivially. 
Thus, we may assume that $K$ is not a complete graph and so by Lemma \ref{maxdeg2tlemma} 
we know that $K$ is not regular and has maximum degree $\Delta=2t$. 
Let $\B = \{B_1, \dots, B_s\}$, let $b_i=|B_i|$ for $i = 1, \dots, s$, and let $M$ be the multiset $M = \{b_1,\dots, b_s\}$. We define $b_{\max} = \max(M)$, $b_{\min} = \min(M)$, and $b_{\min}' = \min(M \setminus \{b_{\min}\})$ (it is possible that $b_{\min}' = b_{\min}$).  
To prove the lemma, it suffices to show that 
\begin{itemize}
\item[(i)] $b_{\max}b_{\min} \geq 2t$; or
\item[(ii)] $b_{\max}b_{\min}' \geq 2t$ and $b_{\min}=|A_x|$ for some $x \in \{1,\dots,r\}$.
\end{itemize}
For then we can take the central vertex of $S$ to 
be $\{B_y\}$ for some $y$ such that $|B_y|=b_{\max}$. 
If (i) holds, then $S$ will be a spanning star, and if (ii) holds, then we can take $V(S)=V(K^\B)\setminus\{A_x\}$ so that 
$V(K^\B)\setminus V(S)=\{A_x\}$. 

Let $|A_i|=a_i$ for $i = 1, \dots, r$ and assume $a_1 \leq \dots \leq a_r$.
For simplicity of notation, let $a = a_1$ be the cardinality of a smallest part,
and let $p=|\{a_i:a_i=a, i=1,\ldots,r\}|$
be the number of parts of cardinality equal to the cardinality of a smallest part.
Observe that $p<r$ because $K$ is not regular.  

If $p=1$, then $a_1<a_2$ and by Lemma \ref{epsboundlemma} we have $a_i=a+1$ for $i=2,\ldots,r$.
Thus, we have $2t=\Delta=b_{\max}=(r-1)(a+1)$ and $b_{\max}b_{\min}\geq 2t$, so (i) holds.  We can henceforth assume $p\geq 2$. 

Our next goal is to prove inequalities (\ref{starlemmaeqn}) and (\ref{starlemmaeqn2}) below, which will be used repeatedly in the remainder of the proof. 
Since a vertex of maximum degree in $K$ is contained in a part of smallest cardinality, 
we have 
\begin{equation}\label{maxdegeqn}
2t = \sum\limits_{i=2}^r a_i.
\end{equation}
Since $p\geq 2$ we have $a_1=a_2$ and so by Lemma \ref{epsboundlemma} we have $a_1+\cdots +a_r\leq r(a+1)-2$. 
Let $\eps_i = 0$ for $i = 1, \dots, p$, and $\eps_i = a_i - (a+1)$ for $i = p+1, \dots, r$.
It follows (using $a_1+\cdots +a_p=ap$ and $a_1+\cdots+a_r\leq r(a+1)-2$) that
\begin{equation}\label{maxeps}
\sum\limits_{i=1}^r \eps_i \leq p-2.
\end{equation}
 
We now show that if 
$I \subseteq \{1,\dots,r\}$ and $\min\{a_i : i \in I\}\geq a+2$, then 
\begin{equation} \label{boundedsumeqn}
\sum\limits_{i\in I} a_i \, \leq \, \frac{\min\{a_i : i \in I\}}{\min\{a_i : i \in I\}-a-1} \sum\limits_{i \in I} \eps_i.
\end{equation}
Write $\sum\limits_{i\in I} a_i$ as 
$c_0a^* + c_1(a^*+1) + \dots + c_q(a^*+q)$ where 
$a^*=\min\{a_i : i \in I\}$, $q$ is given by $a^*+q=\max\{a_i: i\in I\}$, 
and $c_j=|\{a_i: a_i=a^*+j,i\in I\}|$ for $j=0,\ldots,q$,
and define $\eps^*=a^*-(a+1)$.  
Then we have  
$$
\begin{array}{ll}
\sum\limits_{i\in I} a_i&=c_0a^* + c_1(a^*+1) + \dots + c_q(a^*+q)\\
&\leq c_0a^* + c_1a^*(\frac{\eps^*+1}{\eps^*}) + \dots + c_qa^*(\frac{\eps^*+q}{\eps^*})\qquad\text{since } \eps^*\leq a^*\text{ implies }\frac{a^*+j}{a^*}\leq \frac{\eps^*+j}{\eps^*}\\
&=\frac{a^*}{\eps^*}(c_0\eps^* + c_1(\eps^*+1) + \dots + c_q(\eps^*+q))\\
&=\frac{a^*}{\eps^*} \sum\limits_{i \in I} \eps_i
\end{array}
$$
and (\ref{boundedsumeqn}) holds.
It follows easily from (\ref{maxdegeqn}), (\ref{maxeps}) and (\ref{boundedsumeqn}) that if $I \subseteq \{1,\dots,r\}$ and $\min\{a_i : i \in I\}\geq a+2$, then 
\begin{equation} \label{starlemmaeqn}
2t \leq \sum\limits_{i \in \{2,\ldots,r\} \setminus I} a_i + \, \frac{\min\{a_i : i \in I\}}{\min\{a_i : i \in I\}-a-1} (p-2).
\end{equation}
In particular, when $i^* = \min\{i : a_i \geq a+2\}$ and $I = \{i^*, \dots, r\}$, then by (\ref{starlemmaeqn})
we have
\begin{equation} \label{starlemmaeqn2}
2t \leq \sum\limits_{i \in \{2,\ldots,i^*-1\}} a_i + \, (a+2)(p-2)
\end{equation}
because $a_{i^*}\geq a+2$ implies $\frac{a_{i^*}}{a_{i^*}-a-1}\leq a+2$.

The proof now splits into Case 1 where $b_{\max} = ap$ and Case 2 where $b_{\max} > ap$.  
In each case, we show that either (i) or (ii) holds.

\noindent\textbf{Case 1.} Suppose $b_{\max} = ap$.   We consider two subcases, namely $b_{\min} = a+1$ and $b_{\min} \geq a+2$.
\vspace{0.3cm}

\textbf{(a)} Suppose $b_{\min} = a+1$.  First assume $a\geq 2$.
Then there is exactly one part of cardinality $a+1$
and so by (\ref{starlemmaeqn2}) we have
\begin{equation}\label{case1a} 
2t \leq a(p-1) + a+1 + (a+2)(p-2) = 2ap -2a+2p-3. \end{equation}
If $2ap -2a+2p-3 > ap(a+1)$, then $p(a^2-a-2)+2a+3 < 0$, which is a contradiction.
Thus, $2t \leq 2ap -2a+2p-3 \leq ap(a+1) = b_{\max}b_{\min}$, so (i) holds. 

Now assume $a=1$, which implies $b_{\min}=2$ and $b_{\max}=p$. 
We aim to show that (ii) holds with $x = p+1$.  Note that if $p=2$, 
then $b_{\min}=b_{\max}=2$ which implies $K=K_{1^2}$ or $K = K_{1^2,2}$. But $K_{1^2}$ is complete and $K_{1^2,2}$ is not admissible. So $p \geq 3$ and $K$ has exactly one part of cardinality $2$. Thus $b_{\min}' \geq 3$. 
This means that (\ref{case1a}) simplifies to $2t \leq 4p-5<4p$.
Hence we can assume that $b_{\min}' = 3$ because (ii) holds if $b_{\min}' \geq 4$.
If $p=3$, then $b_{\max}=b'_{\min}=3$ and so $K=K_{1^3,2}$ or $K_{1^3,2,3}$. Since neither of these is admissible,
we can assume $p>3$.
Hence $K$ has exactly one part of cardinality $3$, and by (\ref{starlemmaeqn}) with 
$I = \{p+3,\dots,r\}$ we have $2t \leq (p-1) + 2 + 3 + \frac{4}{2}(p-2) = 3p= b_{\max}b'_{\min}$, and again (ii) holds.
\vspace{0.3cm}

\textbf{(b)} Suppose $b_{\min} \geq a+2$.  By (\ref{starlemmaeqn2}), and since $b_{\max} = ap$, we have
\[ 2t \leq a(p-1) + ap + (a+2)(p-2) = 3ap -3a+2p-4. \]
If $3ap -3a+2p-4 > ap(a+2)$, then $p(a^2-a-2)+3a+4 < 0$, which is a contradiction for all $a \geq 2$.  Thus, (i) holds for all $a \geq 2$, and hence we can assume $a=1$.  Now, with $a=1$, we have $b_{\max}=p$ and $2t \leq 5p-7<5p$, and so clearly (i) holds for $b_{\min} \geq 5$. Hence we only need to consider $b_{\min} \in \{3,4\}$.  

If $b_{\min} = 4$, then $K$ has exactly four parts of cardinality $1$, or 
exactly two parts of cardinality $2$, or exactly one part of cardinality $4$.
Moreover, the number of parts having cardinality $3$ is either $0$ or at least $2$. 
If $K$ has exactly four parts of cardinality $1$, then $p=4$ and $b_{\max}=b_{\min}=4$.
This implies that $K=K_{1^4}$, $K=K_{1^4,2^2}$, $K=K_{1^4,4}$, or $K=K_{1^4,2^2,4}$.
The graph $K_{1^4}$ is complete and none of $K_{1^4,2^2}$, $K_{1^4,4}$, or $K_{1^4,2^2,4}$ is admissible.
Thus, we can assume $p\geq 5$. 
If $K$ has exactly two parts of cardinality $2$, then by (\ref{starlemmaeqn2}) we have $2t \leq (p-1) + 4 + 3(p-2) = 4p-3<4p=b_{\max}b_{\min}$ and (i) holds.  
Hence we can assume that $K$ has exactly one part of cardinality $4$ and $b_{\min}' \geq 5$.  Noting that there are at most $b_{\max}=p$ vertices in the parts of cardinality $2$ 
and at most $b_{\max}=p$ vertices in the parts of cardinality $3$, 
by (\ref{starlemmaeqn}) with $i^* = \min \{i : a_i \geq 5\}$ and $I = \{i^*, \dots, r\}$
we have $ 2t \leq (p-1) + p + p + 4 + \frac{5}{3}(p-2) < 5p\leq b_{\max}b'_{\min}$ and so (ii) holds. 

If $b_{\min} = 3$, then $p \geq 3$ (since $a=1$).  Note that if $p=3$ then by (\ref{maxeps}) we have $\sum_{i=1}^r \eps_i \leq 1$ and so $K=K_{1^3}$ or $K=K_{1^3,3}$ 
(recall that $b_{\max} = p =  3$). The graph $K_{1^3}$ is complete and the graph 
$K_{1^3,3}$ is not admissible. Hence we can assume that $p \geq 4$.  
Since $b_{\min}=3$, $K$ has exactly one part of cardinality $3$ and so 
by (\ref{starlemmaeqn}) 
with $i^* = \min \{i : a_i \geq 4\}$
and $I = \{i^*,\dots,r\}$ 
we have $ 2t \leq (p-1) + p + 3 + \frac{4}{2}(p-2) = 4p -2<4p\leq b_{\max}b'_{\min}$.
Thus, (ii) holds.

\vspace{0.3cm}

\noindent \textbf{Case 2.} Suppose $b_{\max} > ap$.  
Let $B_x \in \B$ such that $b_x = b_{\max}$ and let $j$ be such that $a+j$ is the cardinality of the parts of $K$ contained in $B_x$.  Let $\ell = b_x/(a+j)$ so that $\ell$ is the number of parts of cardinality $a+j$ in $K$.

Note that $b_{\min} \geq a+1$ (since $p \geq 2$).   By (\ref{starlemmaeqn2}) we have
\begin{eqnarray}
2t &&\leq  a(p-1) + (a+j)\ell + (a+2)(p-2) \nonumber \\  
   && <  (a+j)\ell \left(3+\frac{2}{a}\right) - 3a - 4,  \label{case2eqn}
\end{eqnarray} 
where the last inequality follows from the fact that $p < (\frac{a+j}{a})\ell$ (recall that $ap < b_{\max} = (a+j)\ell$).
If $(a+j)\ell \left(3+\frac{2}{a}\right) - 3a - 4 > (a+j)\ell(a+1)$, then $(a+j)\ell \left(a-2-\frac{2}{a} \right) + 3a+4 < 0$, which is a contradiction for all $a \geq 3$.  Thus (i) holds for $a \geq 3$, and we only need to consider $a \in \{1,2\}$. 
\vspace{0.3cm}

\textbf{(a)} Suppose $a=2$.  By (\ref{case2eqn}) we have 
$2t \leq 4(2+j)\ell-11\leq 4b_{\max}$ and clearly (i) holds if $b_{\min} \geq 4$.  So we can assume $b_{\min} = 3$, and hence that $K$ has exactly one part of cardinality $3$.  Thus, by (\ref{starlemmaeqn2}), $2t \leq 2(p-1) + 3 + 4(p-2) = 6p - 7<6p<b_{\max}b_{\min}$ (since $b_{\max} > 2p$) and so (i) holds.
\vspace{0.3cm}

\textbf{(b)} Suppose $a=1$.  By (\ref{case2eqn}) we have $2t \leq 5(1+j)\ell-8<5b_{\max}$ and clearly (i) holds if $b_{\min} \geq 5$.  So we can assume $b_{\min} \leq 4$.  In fact, since $p \geq 2$, we have $b_{\min} \in \{2,3,4\}$; we consider these three options in turn.

Suppose $b_{\min}=2$.   First we consider $p \in \{2, 3\}$.  If $p = 2$ then from (\ref{maxeps}) we have $\sum_{i=1}^r \eps_i = 0 $ and hence $K= K_{1^2,2^\ell}$, which is not admissible; if $p=3$ then from (\ref{maxeps}) we have $\sum_{i=1}^r \eps_i \leq 1 $ and hence $K=K_{1^3,2}$ or $K = K_{1^3,2,3}$, which are not admissible.
Hence we can assume that $p \geq 4$.  Since  $b_{\min}=2$, $K$ has exactly one part of cardinality $2$.  Now by (\ref{starlemmaeqn2}) we have $2t \leq p-1+2+3(p-2) = 4p-5<4p<4b_{\max}$.
Clearly (ii) holds if $b_{\min}' \geq 4$. On the other hand, 
if $b_{\min}' = 3$, then $K$ has exactly one part of cardinality $3$ and so by (\ref{starlemmaeqn}) with $i^* = \min \{i : a_i \geq 4\}$ and
$I = \{i^*,\dots,r\}$ we have $2t \leq p-1+2+3+\frac{4}{2}(p-2) = 3p<b_{\max}b'_{\min}$ and again (ii) holds.

Suppose $b_{\min}=3$.   First we consider $p = 3$.  From (\ref{maxeps}) we have $\sum_{i=1}^r \eps_i \leq 1 $ and hence $K=K_{1^3,2^x}$ or $K=K_{1^3,2^x,3}$, for some integer $x$.
But $K_{1^3,2^x,3}$ is is not admissible, so $K=K_{1^3,2^x}$  and (i) holds.  
Hence we can assume that $p \geq 4$.  Since $b_{\min}=3$, $K$ has exactly one part of cardinality $3$, and $b_{\min}' \geq 4$.  Now by (\ref{starlemmaeqn}) with $i^* = \min \{i : a_i \geq 4\}$ and $I = \{i^*,\dots,r\}$, and since $p < b_{\max}$, we have $2t \leq p-1 + b_{\max} + 3 + \frac{4}{2}(p-2) < 4b_{\max}-1$ and thus (ii) holds (since $b'_{\min}\geq 4$).

Suppose $b_{\min}=4$.  First we consider $p=4$.  From (\ref{maxeps}) we have $\sum_{i=1}^r \eps_i \leq 2$, and since $K$ cannot have exactly one part of cardinality $3$ (since $b_{\min}=4$), it follows that $K=K_{1^4,2^x}$ or $K=K_{1^4,2^x,3^2}$ or $K=K_{1^4,2^x,4}$, for some integer $x$, but these are not admissible.
Hence we can assume that $p \geq 5$.  Since $b_{\min}=4$, either $K$ has exactly two parts of cardinality $2$ or exactly one part of cardinality $4$.
If $K$ has exactly two parts of cardinality $2$, then by (\ref{starlemmaeqn2}) we have $2t \leq p-1+4 + 3(p-2) = 4p-3<4p<b_{\max}b_{\min}$ and (i) holds.  So we can assume that $K$ has at least three parts of cardinality $2$ and exactly one part of cardinality $4$ (and so $b_{\min}' \geq 5$).  Now, by (\ref{starlemmaeqn}) 
with $i^*=\min\{i:a_i\geq 3\}$ and $I=\{i^*,\ldots,r\}\setminus\{\ell\}$ 
where $a_\ell=4$, and since $p < b_{\max}$, we have $2t \leq p-1 + b_{\max} + 4 + 3(p-2) < 5b_{\max} - 3$ and (ii) holds (since $b'_{\min}\geq 5$).
\qed

\vspace{0.3cm}

Before proceeding with the next major step in the proof of our main result, namely Lemma \ref{partsizes}, we need the following 
easy consequence of the result of Petersen \cite{Pet} (see \cite{Wes})
that every regular graph of even degree contains a $2$-factor (spanning $2$-regular subgraph), and hence has a $2$-factorisation (decomposition into $2$-factors). 

\begin{lemma}\label{factorisationlemma}
If $K$ is a graph and $t$ is a positive integer such that $2t$ divides $\deg_K(v)$ for each $v\in V(K)$,
then there exists a factorisation $\{F_1,F_2,\ldots,F_t\}$ of $K$ such that $\deg_{F_i}(v)=\frac{\deg_K(v)}t$ for each $v\in V(K)$
and each $i\in\{1,\ldots,t\}$. 
\end{lemma}

\proof
Split each vertex $v$ into $\frac{\deg_K(v)}{2t}$ vertices of degree $2t$ (arbitrarily choosing which edges go with each vertex). The resulting $(2t)$-regular graph has a $2$-factorisation, and recombining the vertices yields the required factorisation of $K$.  
\qed

\begin{lemma}\label{partsizes}
Let $K_{A_1,\ldots,A_r}$ be admissible with $r\geq 2$, let $n=v(K_{A_1,\ldots,A_r})$, let $t=e(K_{A_1,\ldots,A_r})/(n-1)$,
and let $\B$ be the degree partition of $V(K_{A_1,\ldots,A_r})$.
There exists a 
factorisation $\{F_1,\ldots,F_t\}$ of $K_{A_1,\ldots,A_r}$ such that for $i=1,\ldots,t$,
\begin{itemize}
\item [$(1)$] $F_i$ is the vertex-disjoint union of cycles and a single path; and 
\item [$(2)$] $F^\B_i$ is connected.
\end{itemize}
\end{lemma}

\proof
Let $s=|\B|$, let $\B=\{B_1,B_2,\ldots,B_s\}$ and for $j=1,\ldots,s$ let $b_j=|B_j|$.
Let $K^\B$ denote the graph $K^\B_{A_1,\ldots,A_r}$.
Adjoin a new vertex $B_\infty$ to $K^\B$ and for $j=1,\ldots,s$ join $B_\infty$ to the vertex $B_j$ of $K^\B$ with an edge of multiplicity
$2tb_j-\deg_{K^\B}(B_j)$. Let $K^\B_\infty$ be the resulting graph. 
Thus, $\deg_{K^\B_\infty}(B_j)=2tb_j$ for $j=1,\ldots,s$. We now show that  $\deg_{K^\B_\infty}(B_\infty)=2t$. 
By definition, $\deg_{K^\B_\infty}(B_\infty)=\sum_{j=1,\ldots,s}(2tb_j-\deg_{K^\B}(B_j))=2t\sum_{j=1,\ldots,s}b_j-\sum_{j=1,\ldots,s}\deg_{K^\B}(B_j)$. But 
$\sum_{j=1,\ldots,s}\deg_{K^\B}(B_j)=2e(K^\B)=2t(n-1)$ and $\sum_{j=1,\ldots,s}b_j=n$. 
Thus we have $\deg_{K^\B_\infty}(B_\infty)=2tn-2t(n-1)=2t$.
For convenience define $b_\infty=1$ so that $\deg_{K^\B_\infty}(B_j)=2tb_j$ for $j\in\{1,\ldots,s,\infty\}$.
Note that $e(K^\B_\infty)=e(K^\B)+2t=t(n+1)$.

By Lemma \ref{starlemma}, $K^\B$ contains a star $S$ of multiplicity $2t$
such that $V(K^\B)\setminus V(S)$ is either empty or is $\{A_x\}$ for some $x\in\{1,\ldots,r\}$.
Since $2t$ divides the degree of each vertex of $K^\B_\infty$ and each vertex of $S$, 
we have that the degree of each vertex in the graph $K^\B_\infty-E(S)$ 
obtained from $K^\B_\infty$ by deleting the edges of $S$ 
is also divisible by $2t$. Hence by Lemma \ref{factorisationlemma}, $K^\B_\infty-E(S)$ has a factorisation $\{Z_1,\ldots,Z_t\}$ such 
that for $j\in\{1,2,\ldots,s,\infty\}$, $\deg_{Z_x}(B_j)=\deg_{Z_y}(B_j)$ for $1\leq x\leq y\leq t$.
Let $\{S_1,\ldots,S_t\}$ be a factorisation of $S$ where each $S_i$ contains exactly $2$ edges from each set of $2t$ parallel edges
in $S$, and for $i=1,\ldots,t$ let $Z'_i=Z_i\cup S_i$ so that $\{Z'_1,\ldots,Z'_t\}$ is a factorisation of $K^\B_\infty$. 
It follows that for $j\in\{1,2,\ldots,s,\infty\}$, $\deg_{Z'_x}(B_j)=\deg_{Z'_y}(B_j)$ for $1\leq x\leq y\leq t$,
and hence that for $i=1,\dots,t$ and $j\in\{1,2,\ldots,s,\infty\}$, $\deg_{Z'_i}(B_j)=2b_j$. 
Thus, for $i=1,\dots,t$ we have $e(Z'_i)=e(K^\B_\infty)/t=n+1$.

For $i=1,\ldots,t$ let $H_i$ be the factor of $K^\B$ obtained by deleting $B_\infty$ from $Z'_i$. 
Since each $H_i$ contains $S_i$ as a subgraph, 
we have $e(Z'_i)=n+1$, $\deg_{Z'_i}(B_j)=2b_j$ and $\deg_{Z'_i}(B_\infty)=2$,
and it follows that for $i=1,\ldots,t$, 
\begin{itemize}
\item [(a)] $\deg_{H_i}(B_j)\leq 2b_j$ for $j=1\ldots,s$;
\item [(b)] $e(H_i)=n-1$; and
\item [(c)] $H_i$ is connected.
\end{itemize}
To see that each $H_i$ is connected, observe that $H_i$ contains the edges of $S_i$. Thus it is clear that $H_i$ is
connected if $S_i$ is a spanning star.
If $S_i$ is not a spanning star, then there is a component of
$H_i$ (namely, the component of $H_i$ containing the edges of $S_i$)
containing all the vertices of $K^\B$, except possibly some vertex $B_x=A_x$.
Thus, $H_i$ is again connected unless it is the case that in $Z'_i$, all the edges incident with $A_x$ are joined to $B_\infty$.
This could happen only if $Z'_i$ has a $2$-cycle on $A_x$ and $B_\infty$ 
(because $\deg_{Z'_i}(B_\infty)=2$)
and $|A_x|=1$  
(because we have $A_x=B_x$ and $\deg_{Z'_i}(B_x)=2b_x$ which implies $b_x=1$). But if $|A_x|=1$, then $A_x$ is not joined to
$B_\infty$ in $K^\B_\infty$. So this situation does not arise.

Now, for $i=1,\ldots,t$, let $G_i$ be the factor of $K_{A_1,\ldots,A_r}$ 
that contains the edge $e$ of $K_{A_1,\ldots,A_r}$ if and only if $e^\B$ is
in $H_i$.
Thus, for $i=1,\ldots,t$, we have $G^\B_i=H_i$ and it 
follows that $\{G_1,\ldots,G_t\}$ is a factorisation of $K_{A_1,\ldots,A_r}$ such that for $i=1,\ldots,t$, 
\begin{itemize}
\item [(d)] $\sum_{v\in B_j}\deg_{G_i}(v)\leq 2b_j$ for $j=1,\ldots,s$; 
\item [(e)] $e(G_i)=n-1$; and
\item [(f)] $G^\B_i$ is connected.
\end{itemize}

Now let $\{F_1,\ldots,F_t\}$ be the factorisation of $K_{A_1,\ldots,A_r}$ obtained from $\{G_1,\ldots,G_t\}$ by applying 
Lemma \ref{equitabilise}. 
It follows from (d) and (and from (2) in Lemma \ref{equitabilise}) that for $i=1,\ldots,t$, we have 
$\deg_{F_i}(v)\leq 2$ for each vertex $v$ of $K_{A_1,\ldots,A_r}$. Thus, since $e(F_i)=e(G_i)=n-1$, each $F_i$ is the vertex-disjoint union
of cycles and a single path. 
Finally, it follows from (f) (and from (1) in Lemma \ref{equitabilise}) that $F^\B_i$ is connected.
Thus, $\{F_1,\ldots,F_t\}$ is the required factorisation of $K_{A_1,\ldots,A_r}$.
\qed

\begin{lemma}\label{partsgraph}
Let $K_{A_1,\ldots,A_r}$ be admissible with $r\geq 2$, let $n=v(K_{A_1,\ldots,A_r})$, let $t=e(K_{A_1,\ldots,A_r})/(n-1)$ and 
let $\A=\{A_1,\ldots,A_r\}$.
There exists a 
factorisation $\{F_1,\ldots,F_t\}$ of $K_{A_1,\ldots,A_r}$ such that for $i=1,\ldots,t$,
\begin{itemize}
\item [$(1)$] $F_i$ is the vertex-disjoint union of cycles and a single path; and
\item [$(2)$] $F^\A_i$ is connected.
\end{itemize}
\end{lemma}

\proof
Let $\B$ be the degree partition of $K_{A_1,\ldots,A_r}$ and let 
$\{H_1,\ldots,H_t\}$ be a factorisation $K_{A_1,\ldots,A_r}$ such that for $i=1,\ldots,t$, 
$H_i$ is the vertex-disjoint union of cycles and a single path and $H^\B_i$ is connected. 
Such a factorisation exists by Lemma \ref{partsizes}.
Let $C=\sum_{i=1,\ldots,t}c(H^\A_i)$ so that $C$ is the number of components in the factors of 
the factorisation $\{H^\A_1,\ldots,H^\A_t\}$.
If $C=t$, then each $H^\A_i$ is connected and letting $F_i=H_i$ for $i=1,\ldots,t$ gives the required factorisation. 
Thus, we can assume $C>t$, which means that there exists a $k$ such that $H^\A_k$ is disconnected. 

Since $H^\B_k$ is connected this implies that there exist $x,y\in\{1,\ldots,r\}$ such that  
$|A_x|=|A_y|$ and $A_x$ and $A_y$ are in distinct components of $H^\A_k$. 
We now show that we can apply Lemma \ref{connectinglemma} with $\a=A_x$ and $\b=A_y$
to the factorisation $\{H^\A_1,\ldots,H^\A_t\}$ of $K^\A_{A_1,\ldots,A_r}$.
Clearly, the transposition $(A_x\,A_y)$ is an automorphism of $K^\A_{A_1,\ldots,A_r}$,
and it follows from the fact that $H_i$ is the vertex-disjoint union of cycles and a single path that 
for $i=1,\ldots,t$ either 
$A_x$ and $A_y$ are in the same component of $H^\A_i$ or 
at least one of $A_x$ and $A_y$ is in a cycle of $H^\A_i$.
Thus, we can indeed apply Lemma \ref{connectinglemma}, and we let $\{G_1,\ldots,G_t\}$ be the resulting factorisation of 
$K^\A_{A_1,\ldots,A_r}$.

Let $C'=\sum_{i=1,\ldots,t}c(G_i)$.
It follows from (1)-(3) and (5) of Lemma \ref{connectinglemma} that 
$c(G_i)\leq c(H^\A_i)$ for $i=1,\ldots,t$ and $c(G_k)<c(H^\A_k)$.
Thus, $C'<C$. 
For $i=1,\ldots,t$, let $X_i$ be the factor of $K_{A_1,\ldots,A_r}$ that contains the edge
$e$ of $K_{A_1,\ldots,A_r}$ if and only if $e^\A$ is
in $G_i$, and then let $\{Y_1,\ldots,Y_t\}$ be the factorisation of $K_{A_1,\ldots,A_r}$ obtained from
$\{X_1,\ldots,X_t\}$ by applying Lemma \ref{equitabilise}.
Thus, for $i=1,\ldots,t$ we have $Y^\A_i=X^\A_i=G_i$ and $Y_i$ is almost regular on each $A_j$. 
Note that $\sum_{i=1,\ldots,t}c(Y^\A_i)=C'<C$ (since $Y^\A_i=G_i$). 

We will show that each $Y_i$ is the vertex-disjoint union of cycles and a single path, 
and that each $Y^\B_i$ is connected. 
Thus, if $C'=t$, then $\{Y_1,\ldots,Y_t\}$ is the required factorisation of $K_{A_1,\ldots,A_r}$.
Otherwise, we can repeat the above procedure (each time reducing the number of components in the factors of the factorisation) 
until we obtain the required factorisation $\{F_1,\ldots,F_t\}$ 
of $K_{A_1,\ldots,A_r}$.

For $i=1,\ldots,t$ and $j=1,\ldots,r$, we have $\sum_{v\in A_j}\deg_{Y_i}(v)=\deg_{Y^\A_i}(A_j)=\deg_{G_i}(A_j)$.
Thus, since each $Y_i$ is almost regular on each $A_j$, if we can show that 
$\deg_{G_i}(A_j)\leq 2|A_j|$
then we have that each $Y_i$ has maximum degree at most $2$. 
Since we know that $e(Y_i)=e(H_i)=n-1$, this implies that each $Y_i$ is the vertex-disjoint union of cycles and a single path.

Now, for $i=1,\ldots,t$ and $j=1,\ldots,r$, we have $\deg_{H^\A_i}(A_j)=\sum_{v\in A_j}\deg_{H_i}(v)\leq 2|A_j|$
(because $H_i$ is the vertex-disjoint union of cycles and a single path).
But for $j\notin\{x,y\}$ and for $i=1,\ldots,t$, we have $\deg_{G_i}(A_j)=\deg_{H^\A_i}(A_j)$ and so we have 
$\deg_{G_i}(A_j)\leq 2|A_j|$ as required for $j\notin\{x,y\}$.
For $i=1,\ldots,t$ we also have $\deg_{G_i}(A_x)+\deg_{G_i}(A_y)=\deg_{H^\A_i}(A_x)+\deg_{H^\A_i}(A_y)\leq 2|A_x|+2|A_y|$
(the equality $\deg_{G_i}(A_x)+\deg_{G_i}(A_y)=\deg_{H^\A_i}(A_x)+\deg_{H^\A_i}(A_y)$
following from (1)-(3) of Lemma \ref{connectinglemma}).
It follows from $\deg_{G_i}(A_x)+\deg_{G_i}(A_y)\leq 2|A_x|+2|A_y|$ and the fact that each of $G_1,\ldots,G_t$ is almost regular on 
$\{A_x,A_y\}$ that $\deg_{G_i}(A_x)\leq 2|A_x|$ and $\deg_{G_i}(A_y)\leq 2|A_y|$ (recall that $|A_x|=|A_y|$).
Thus, each $Y_i$ is indeed the vertex-disjoint union of cycles and a single path.

It remains only to show that each $Y^\B_i$ is connected. 
Since $Y^\A_i=X^\A_i=G_i$ and
$G_i-\{A_x,A_y\}=H^\A_i-\{A_x,A_y\}$ (by (3) of Lemma \ref{connectinglemma}),
we have $Y^\A_i-\{A_x,A_y\}=H^\A_i-\{A_x,A_y\}$,
and it follows that $Y^\B_i=H^\B_i$ (because the elements of $A_x\cup A_y$ are all in the same part of $\B$). Thus, $Y^\B_i$ is connected because $H^\B_i$ is connected.
\qed

\begin{theorem}\label{mainresult}
A complete multipartite graph $K$ with $n>1$ vertices and $m$ edges 
can be decomposed into edge-disjoint Hamilton paths if and only if 
$\frac m{n-1}$ is an integer and the maximum degree of $K$ is at most $\frac {2m}{n-1}$.
\end{theorem}

\proof
It was noted earlier in the paper that if $K$ has a Hamilton path decomposition, then $t=\frac m{n-1}$ is an integer and $K$ has maximum degree at most $2t$. We now show that these conditions are sufficient for the existence of a Hamilton path decomposition of $K$. The proof is somewhat similar to the proof of Lemma \ref{partsgraph}.
Let $K\cong K_{A_1,\ldots,A_r}$ be admissible, and let $\{H_1,\ldots,H_t\}$ be a factorisation of $K_{A_1,\ldots,A_r}$ such that each 
$H_i$ is the vertex-disjoint union of cycles and a single path and each $H^\A_i$ is connected. 
Such a factorisation exists by Lemma \ref{partsgraph}. 

Let $C=\sum_{i=1,\ldots,t}c(H_i)$ so that $C$ is the number of components in the factors of 
the factorisation $\{H_1,\ldots,H_t\}$.
If $C=t$, then each $H_i$ is connected and letting $F_i=H_i$ for $i=1,\ldots,t$ gives the required Hamilton path decomposition. 
Thus, we can assume $C>t$, which means that there exists a $k$ such that $H_k$ is disconnected. 
Since $H^\A_k$ is connected this implies that there exist $x\in \{1,\ldots,r\}$ and $u,v\in A_x$ such that  
$u$ and $v$ are in distinct components of $H_k$. 

We now show that we can apply Lemma \ref{connectinglemma} with $\a=u$ and $\b=v$
to the factorisation $\{H_1,\ldots,H_t\}$ of $K_{A_1,\ldots,A_r}$.
Clearly, the transposition $(u\,v)$ is an automorphism of $K_{A_1,\ldots,A_r}$,
and since $H_i$ is the vertex-disjoint union of cycles and a single path, for $i=1,\ldots,t$ we have either 
$u$ and $v$ are in the same component of $H_i$ or 
at least one of $u$ and $v$ is in a cycle of $H_i$.
Thus, we can indeed apply Lemma \ref{connectinglemma}, and we let $\{G_1,\ldots,G_t\}$ be the resulting factorisation of 
$K_{A_1,\ldots,A_r}$. Observe that $G^\A_i=H^A_i$ for $i=1,\ldots,t$.

It is clear that each $G_i$ is the vertex-disjoint union of cycles and a single path, 
and that each $G^\A_i$ is connected (because $G^\A_i=H^A_i$). 
Let $C'=\sum_{i=1,\ldots,t}c(G_i)$.
It follows from (1)-(3) and (5) of Lemma \ref{connectinglemma} that 
$c(G_i)\leq c(H_i)$ for $i=1,\ldots,t$ and $c(G_k)<c(H_k)$.
Thus, $C'<C$. 
If $C'=t$, then $\{G_1,\ldots,G_t\}$ is the required Hamilton path decomposition of $K_{A_1,\ldots,A_r}$.
Otherwise, we can repeat the above procedure (each time reducing the number of components in the factors of the factorisation) 
until we obtain the required Hamilton path decomposition of $K_{A_1,\ldots,A_r}$.
\qed

  \let\oldthebibliography=\thebibliography
  \let\endoldthebibliography=\endthebibliography
  \renewenvironment{thebibliography}[1]{%
    \begin{oldthebibliography}{#1}%
      \setlength{\parskip}{0.4ex plus 0.1ex minus 0.1ex}%
      \setlength{\itemsep}{0.4ex plus 0.1ex minus 0.1ex}%
  }%
  {%
    \end{oldthebibliography}%
  }


\begin{thebibliography}{99}
\def\baselinestretch{1}\small\normalsize
\begin{small}



\bibitem{Als}
B. Alspach, 
The wonderful Walecki construction, 
{\it Bull. Inst. Combin. Appl.}
{\bf 52} (2008) 7--20.

\bibitem{AndRod}
L. D. Andersen and C. A. Rodger,
Decompositions of complete graphs: embedding partial edge-colourings and
the method of amalgamations,
Surveys in Combinatorics, London Mathematical Society Lecture Note Series 
{\bf 307} (2003), 7--41.

\bibitem{BotJim}
F. Botler and A. Jim\'{e}nez,
On path decompositions of $2k$-regular graphs,
{\it Discrete Math.}
{\bf 340} (2017) 1405--1411.


\bibitem{BotMotOshWak}
F. Botler, G. O. Mota, M. T. I. Oshiro, Y. Wakabayashi, 
Decomposing highly edge-connected graphs into paths of any given length,
{\it J. Combin. Theory Ser. B}
{\bf 122} (2017), 508--542. 



\bibitem{Bry}
D. Bryant, 
Packing paths in complete graphs,
{\it J. Combin. Theory Ser. B}
{\bf 100} (2010) 206--215.

\bibitem{Bry2}
D. Bryant,
On Almost-Regular Edge Colourings of Hypergraphs, 
{\it Electron. J. Combin.} 
{\bf 23 (4)} (2016) \#P4.7 7 pp.



\bibitem{BryMae} 
D. Bryant and B. Maenhaut, 
Almost regular edge colourings and regular decompositions of complete graphs, 
{\it J. Combin. Des.}, 
{\bf 16} (2008) 499--506.

    
\bibitem{HarMcG}
P. Harding and S. McGuinness,
Gallai's Conjecture for graphs of girth at least four,
{\it J. Graph Theory}
{\bf 75} (2014), no. 3, 256--274. 

\bibitem{HofRod} 
D. G. Hoffman and C.A.~Rodger, 
The chromatic index of complete multipartite graphs, 
{\it J. Graph Theory}, {\bf 16} (1992), 159--163.

\bibitem{LasAue}
R. Laskar and B. Auerbach, On decomopositions of $r$-partite graphs into edge-disjoint Hamilton circuits,
{\it Discrete Math.},
{\bf 14} (1976) 265--268.

\bibitem{Lov}
L. Lov\'{a}sz, On covering of graphs, in: P. Erd\"{o}s, G. Katona (Eds.), Theory of Graphs, Academic
Press, New York, 1968, pp. 231--236.


\bibitem{Luc}
E. Lucas,
R\'{e}cr\'{e}ations Math\'{e}matiques,
Volume II, Gauthiers-Villars, Paris, 1894.





\bibitem{Par}
C. A. Parker, 
Complete bipartite graph path decompositions, 
Ph.D. Thesis, Auburn University, 1998.

\bibitem{Pet}
J. Petersen, 
Die Theorie der regul\"aren graphs, {\it Acta Math.}, {\bf 15} (1891), 193--220. 

\bibitem{Tar}
M. Tarsi,
Decompositions of a complete multipartite graph into simple paths: Non balanced handcuffed designs,
{\it J. Combin. Theory Ser. A}
{\bf 34} (1983), 60--70.

\bibitem{Tho1}
C. Thomassen, 
Decompositions of highly connected graphs into paths of length 3, 
{\it J. Graph Theory}
{\bf 58} (2008) 286--292.

\bibitem{Tho2}
C. Thomassen, 
Edge-decompositions of highly connected graphs into paths,
{\it Abh. Math. Semin. Univ. Hambg.}
{\bf 78} (2008), no. 1, 17--26. 

\bibitem{Wes}
D. B. West, 
Introduction to Graph Theory - Second edition, Published by Prentice Hall 1996, 2001.







\end{small}
\end{thebibliography}
\end{document}